# To the Many-Hilbert-Space Theory of Quantum Measurements

by

## Z.I.ISMAILOV and E.OTKUN ÇEVIK

Karadeniz Technical University, Faculty of Sciences, Department of Mathematics
61080 Trabzon, TURKEY
e-mail address : zameddin@yahoo.com

**Abstract:** In this work, a connection between some spectral properties of direct integral of operators in the direct integral of Hilbert spaces and their coordinate operators has been investigated.

## 1. Introduction

It is known that the general theory of linear closed operators in Hilbert spaces and its applications to physical problems have been investigated by many researchers(for example,see [1,2]).But many physical problems require studying the theory of linear operators in direct sums or in general direct integrals of Hilbert spaces. The concepts of direct integral of Hilbert spaces and direct integral of operators as a generalization of the concept of direct sum of Hilbert spaces and direct sum of operators were introduced to mathematics and developed in 1949 by John von Neumann in his work [3] .These subjects were incorporated in several works(see [4 − 7]).A spectral theory of some operators on a finite sum of Hilbert spaces was investigated by N.Dunford [8,9].Note that,in terms of application ,there are some results in papers [10 − 12] in the finite sum cases .Also in the infinite direct sum cases there is a work[13] ,in which some spectral and compactness properties are surveyed.

Furthermore,some spectral investigations of the direct integral of operators in the direct integral of Hilbert spaces have been provided by T.R.Chow[14],T.R.Chow,F.Gilfeather[15], E.A.Azoff[16,17]L.A.Fialkow[18].

It must be noted that the theory of direct integral of Hilbert spaces and operators on the these spaces has important role in the representation theory of locally compact groups,in the theory of decomposition rings of operators to factors,invariant measures,reduction theory,von Neumann algebras and ets.On the other hand, many physical problems of today arising in the modelling of processes of multiparticle quantum mechanics , quantum field theory and in the physics of rigid bodies support to study a theory of  direct integral of operators in the direct integral of Hilbert spaces (see [19] and references in it).

Numereous scientific investigations have been done to  explain quantum measurements .Dealing with these subjects, S.Machida and M.Namiki [20 − 22] (see, also  [23]) have offered many-Hilbert-space theory lately.Also,they find that decoherence of the wave function is only necessity to formulate quantum mechanical measurements. A direct integral space of continiously many Hilbert spaces and a continious superselection rule is starring in their theory.The direct integral space structure is assigned to the measurement devise which reflects its macroscopic features.On the other hand,observed system is left to be defined by a single Hilbert space.Furthermore,they have examined double slit experiments and negative result experiments in the framework of the many-Hilbert-space theory.Note that a direct integral space of continiously many Hilbert spaces often arises in the quantum version of Lax-Phillips theory[24].In this investigation the direct integral space is introduced in order to the allow the generator of motion to have a spectrum in real axis ,which is a necessary condition for the application of the Lax-Phillips theory.

In second section of this paper a connection between parts of spectrum, resolvent set  of direct integral of operators defined in the direct  integral of Hilbert spaces  and   parts of the spectrum of ''coordinate operators ''has been established . Note that the another approach to analogous problem has been  used in work [14] . In this present paper sharp formulas for the connection are given. In third section the compactness properties of these operators have been  examined.Finally,in fourth section in special case the analogous questions  for the power and polynomially bounded operators have been provided.

In the special case of direct sum of Hilbert spaces,  these questions have been investigated in [ 13].

Along this paper   the triplet   $(\Lambda, \Sigma, \mu)$   be any measure space and  the Hilbert spaces are looked at will become infinite dimensional. In addition ,the space of compact operators and Schatten-von Neumann classes in any Hilbert space will be denoted by $C_\infty(.)$ and $C_p(.), 1 \leq p < \infty$, respectively.

## 2. On the spectrum of  direct  integral of operators

In this section, the relationship between  the spectrum and resolvent sets of the direct integral of operators and its coordinate operators will be investigated.

Before of all prove the following result.

**Theorem 2.1**.For the operator $A = (A_\lambda), \lambda \in \Lambda, A = \int_\Lambda^\oplus A_\lambda \, d\mu(\lambda)$ in the Hilbert space $H=(H_\lambda), \lambda \in \Lambda$, $H = \int_\Lambda^\oplus H_\lambda \, d\mu(\lambda)$   are true

$$\sigma_p(A) \subset \bigcup_{\lambda \in \Lambda} \sigma_p(A_\lambda),$$

$$\{\tau \in \bigcap_{\lambda \in \Lambda} \sigma_p(A_\lambda): A_\lambda x_\lambda^\tau = \tau x_\lambda^\tau, x_\lambda^\tau \in D(A_\lambda), x_\lambda^\tau \neq 0, \|x_\lambda^\tau\| \in L^2(\Lambda)\} \subset \sigma_p(A)$$

**Proof**. For any $\tau \in \sigma_p(A)$ there exist element $x = (x_\lambda), \lambda \in \Lambda$ such that $x \neq 0, x \in D(A)$ and $A(x_\lambda) = \tau(x_\lambda)$. Then almost everywhere $\lambda \in \Lambda$ with respect to measure µ it is true

$$A_\lambda x_\lambda = \tau x_\lambda$$

Since $x \neq 0$, then there exist $\lambda_* \in \Lambda$ which satisfy the above equality and $x_{\lambda_*} \in D(A_{\lambda_*}), x_{\lambda_*} \neq 0$.

This means that $\tau \in \sigma_p(A_{\lambda_*})$. Hence $\tau \in \bigcup_{\lambda \in \Lambda} \sigma_p(A_\lambda)$ and from this it is obtained that

$$\sigma_p(A) \subset \bigcup_{\lambda \in \Lambda} \sigma_p(A_\lambda)$$

The proof of the second proposition is clear.

Actually, in one special case the following stronger assertions are true.

**Theorem 2.3.** Assume that every one-point set is measurable and its measure is positive. For the parts of spectrum and resolvent sets of the operator $A = (A_\lambda), \lambda \in \Lambda, A = \int_\Lambda^\oplus A_\lambda \, d\mu(\lambda), A: D(A) \subset H \to H$ in Hilbert space $H = (H_\lambda), \lambda \in \Lambda, H = \int_\Lambda^\oplus H_\lambda \, d\mu(\lambda)$ the following claims are true

(1) $\sigma_p(A) = \bigcup_{\lambda \in \Lambda} \sigma_p(A_\lambda)$;

(2) $\sigma_c(A) = \{(\bigcap_{\lambda \in \Lambda}[\sigma_c(A_\lambda) \cup \rho(A_\lambda)]) \cap (\bigcup_{\lambda \in \Lambda} \sigma_c(A_\lambda))\} \cup \{\tau \in \bigcap_{\lambda \in \Lambda} \rho(A_\lambda): \sup \|R_\tau(A_\lambda)\| = \infty\}$;

(3) $\sigma_r(A) = (\bigcap_{\lambda \in \Lambda}[\sigma_c(A_\lambda) \cup \sigma_r(A_\lambda) \cup \rho(A_\lambda)]) \cap (\bigcup_{\lambda \in \Lambda} \sigma_r(A_\lambda))$;

(4) $\rho(A) = \{\tau \in \bigcap_{\lambda \in \Lambda} \rho(A_\lambda): \sup \|R_\tau(A_\lambda)\| < \infty\}$;

**Proof.** Here only the first and second relations of theorem will be prove. The validity of other propositions can be proved by the similar ideas.

Assumed that $\tau \in \sigma_p(A)$. Then there exist $x = (x_\lambda^\tau) \neq 0$, $(x_\lambda^\tau) \in D(A)$ such that $Ax = \tau x$. So for every $\lambda \in \Lambda$, $A_\lambda x_\lambda^\tau = \tau x_\lambda^\tau, x_\lambda^\tau \in D(A_\lambda)$ and for some $\lambda^* \in \Lambda, x_{\lambda^*}^\tau \neq 0$. Hence $\tau \in \sigma_p(A_{\lambda^*})$ and from this $\tau \in \bigcup_{\lambda \in \Lambda} \sigma_p(A_\lambda)$. On the contrary, assumed that $\tau \in \bigcup_{\lambda \in \Lambda} \sigma_p(A_\lambda)$. Then at least one index $\lambda^* \in \Lambda$ is hold $\tau \in \sigma_p(A_{\lambda^*})$, i.e. for some $x_{\lambda^*}^\tau \neq 0, x_{\lambda^*}^\tau \in D(A_{\lambda^*})$ is true $A_{\lambda^*} x_{\lambda^*}^\tau = \tau x_{\lambda^*}^\tau$. In this case for the element $x = (x_\lambda) \neq 0$, $(x_\lambda) \in D(A), \lambda \neq \lambda^*$, $x_\lambda = 0$ and $x_{\lambda^*} = x_{\lambda^*}^\tau$ we have $Ax = \tau x$. Consequently, $\tau \in \sigma_p(A)$.

Now we prove the second relation on the continuous spectrum. Let $\tau \in \sigma_c(A)$. In this case by the definition of continuous spectrum $A - \tau E$ is a one-to-one operator, $R(A - \tau E) \neq H$ and $R(A - \tau E)$ is dense in $H$. From this and definition of direct integral it implies that for every $\lambda \in \Lambda$ an operator $A_\lambda - \tau E_\lambda$ is a one-to-one operator in $H_\lambda$, $R(A_\lambda - \tau E_\lambda) \neq H_\lambda$ and $R(A_\lambda - \tau E\lambda)$ is dense in $H\lambda$. Since $\tau \in \sigma c A$, then

$\tau \in (\bigcap_{\lambda \in \Lambda} (\sigma_c(A_\lambda) \cup \rho(A_\lambda))) \cap (\bigcup_{\lambda \in \Lambda} \sigma_c(A_\lambda))$ or $\tau \in \bigcap_{\lambda \in \Lambda \backslash \Omega_\tau} \rho(A_\lambda): \sup \|R_\tau(A_\lambda)\| = \infty$.

This means that

$$\sigma_c(A) \subset [(\cap_{\lambda \in \Lambda}(\sigma_c(A_\lambda) \cup \rho(A_\lambda))) \cap (\cup_{\lambda \in \Lambda} \sigma_c(A_\lambda))] \cup \{\tau \in \cap_{\lambda \in \Lambda} \rho(A_\lambda): \sup\|R_\tau(A_\lambda)\| = \infty\}.$$

It is easy to prove the inverse implication.

On the other hand the simple calculations show that the following relations are true.

**Corollary 2.4**. Under the assumptions of last theorem we have

$$\sigma_c(A) = (\cup_{\lambda \in \Lambda} \sigma_p(A_\lambda))^c \cap (\cup_{\lambda \in \Lambda} \sigma_r(A_\lambda))^c \cap (\cup_{\lambda \in \Lambda} \sigma_c(A_\lambda)) \cup \{\tau \in \cap_{\lambda \in \Lambda} \rho(A_\lambda): \sup\|R_\tau(A_\lambda)\| = \infty \;,$$

$$\sigma_r(A) = (\cup_{\lambda \in \Lambda} \sigma_p(A_\lambda))^c \cap (\cup_{\lambda \in \Lambda} \sigma_r(A_\lambda)) \;.$$

**Corollary 2.5.** Let $\Lambda = \{\lambda_1, \lambda_2, \ldots\ldots\ldots\lambda_n\}, n \leq \infty$ be any countable set, $\Sigma = P(\Lambda)$ and $\mu$ be any measure with property $\mu(\{\lambda\}) > 0$ for every point $\lambda \in \Lambda$. In this case the formulas

$$\sigma_p(A) = \cup_{m=1}^n \sigma_p(A_{\lambda_m}),$$

$$\sigma_c(A) = [(\cap_{m=1}^n(\sigma_c(A_{\lambda_m}) \cup \rho(A_{\lambda_m}))) \cap (\cup_{m=1}^n \sigma_c(A_{\lambda_m}))] \cup [\cap_{m=1}^n \{\tau \in \rho(A_{\lambda_m}): \sup\|R_\tau(A_{\lambda_m})\| = \infty\}],$$

$$\sigma_r(A) = (\cap_{m=1}^n(\sigma_c(A_{\lambda_m}) \cup \sigma_r(A_{\lambda_m}) \cup \rho(A_{\lambda_m}))) \cap (\cup_{m=1}^n \sigma_r(A_{\lambda_m})),$$

$$\rho(A) = \{\tau \in \cap_{m=1}^n \rho(A_{\lambda_m}): \sup\{\|R_\tau(A_{\lambda_m})\|: 1 \leq m \leq n\} < \infty\}.$$

are true.

Note that when $\Lambda = N, \Sigma = P(N)$, $\mu$ is counting measure the analogous results have been established in work [13].

### 3. SOME COMPACTNESS PROPERTIES of DIRECT INTEGRAL of OPERATORS

In this section the compactness and spectral properties between direct integral of operators and their "coordinate operators" have been established. In general, there is not any relation between mentioned operators in compactness means.

**Example 3.1**. $\Lambda = N, \Sigma = P(N), \mu$ -counting measure, $H_n = \mathbb{C}, A_n = \mathbb{C}, n \geq 1, H = \oplus_{n=1}^\infty H_n, A = \oplus_{n=1}^\infty A_n$.

In this case for every $n \geq 1, A_n \in C_\infty(H_n)$, but $A \not\in C_\infty(H)$.

**Example 3.2**. In some cases from the relations $A \in C_\infty(H)$ no implies $A_n \in C_\infty(H_n)$ for every $n \geq 1$.

Indeed, from the definition of direct integral of operators on the set having null $\mu$-measure the coordinate operators may be defined by arbitrary way.

But in certain situations there are concrete results.

**Theorem 3.3.** Let $\Lambda = \{\lambda_1, \lambda_2, \ldots \ldots \ldots \lambda_n\}, n \leq \infty$ be any countable set, $\Sigma = P(\Lambda)$ and $\mu$ be any measure with property $\mu(\{\lambda\}) > 0$ for every point $\lambda \in \Lambda$. Then

(1) If $A = \oplus_{m=1}^n A_{\lambda_m} \in C_\infty(H), H = \oplus_{m=1}^n H_{\lambda_m}$, then for every $1 \leq m \leq n$, $A_{\lambda_m} \in C_\infty(H_{\lambda_m})$.

(2) Let $\Lambda$ infinite countable set and for every $n \geq 1, A_{\lambda_n} \in C_\infty(H_{\lambda_n})$. In this case

$$A = \oplus_{n=1}^\infty A_{\lambda_n} \in C_\infty(H) \text{ if and only if } \lim_{n\to\infty} \|A_{\lambda_n}\| = 0.$$

This theorem is proved by analogous scheme of the proof in theorem 4.6 in [13].

Now give one characterizating theorem on the point spectrum of compact direct integral of operators.

**Theorem 3.4.** Assumed that $A = (A_\lambda), \lambda \in \Lambda, A = \int_\Lambda^\oplus A_\lambda \, d\mu(\lambda)$ in the Hilbert space $H=(H_\lambda), \lambda \in \Lambda$, $H = \int_\Lambda^\oplus H_\lambda \, d\mu(\lambda)$, $A_\lambda \in C_\infty(H_\lambda), \lambda \in \Lambda$ and $A \in C_\infty(H)$. In this case there exist countable subset $\Lambda^* = \{\lambda_1, \lambda_2, \lambda_3, \ldots, \lambda_n\}, n \leq \infty$ of $\Lambda$ such that the set $\Lambda^*$ is minimal and

$$\sigma_p(A) = \cup_{m=1}^n \sigma_p(A_{\lambda_m}).$$

From the definition of singular number $s(.)$ (or characteristic numbers) of any compact operator in any Hilbert space [1] and Theorems 2.1 and 3.4 it is easy to prove the validity of the following result.

**Theorem 3.5.** Assumed that $A = (A_\lambda), \lambda \in \Lambda, A = \int_\Lambda^\oplus A_\lambda \, d\mu(\lambda)$ in the Hilbert space $H=(H_\lambda), \lambda \in \Lambda$, $H = \int_\Lambda^\oplus H_\lambda \, d\mu(\lambda)$, $A_\lambda \in C_\infty(H_\lambda), \lambda \in \Lambda$ and $A \in C_\infty(H)$. In this case there exist countable subset $\Lambda^* = \{\lambda_1, \lambda_2, \lambda_3, \ldots, \lambda_n\}, n \leq \infty$ of $\Lambda$ such that

(1) $\{s_k(A): k \geq 1\} = \cup_{m=1}^n \{s_q(A_{\lambda_m}): q \geq 1\}$;

(2) if $A \in C_p(H), 1 \leq p < \infty$, then for every $m, 1 \leq m \leq n$, $A_{\lambda_m} \in C_p(H_{\lambda_m})$;

(3) Let $A_{\lambda_m} \in C_{p(\lambda_m)}(H_{\lambda_m}), 1 \leq m \leq n, 1 \leq p(\lambda_m) < \infty, 1 \leq p = sup\{p(\lambda_m): 1 \leq m \leq n < \infty$. Then $A \in C_p H$ if and only if the series $m=1nq=1\infty sqp(A\lambda m)$ converges.

(4) If $A_{\lambda_m} \in C_{p(\lambda_m)}(H_{\lambda_m}), 1 \leq m \leq n, 1 \leq p(\lambda_m) < \infty, p = sup\{p(\lambda_m): 1 \leq m \leq n\} < \infty$ and the series $\sum_{m=1}^n \sum_{q=1}^\infty s_q^{p(\lambda_m)}(A_{\lambda_m})$ is convergent, then $A \in C_p(H)$.

(5) If $A_{\lambda_m} \in C_{p(\lambda_m)}(H_{\lambda_m}), 1 \leq m \leq n, 1 \leq p(\lambda_m) < \infty, p(\lambda_m) = inf\{\alpha \in [1, \infty): A_{\lambda_m} \in C\alpha H\lambda m$ and $supp\lambda m: 1\leq m\leq n=\infty$, then for every $1\leq p <\infty$, $A\bar\in Cp(H)$.

(6) If $A_{\lambda_m} \in C_{p(\lambda_m)}(H_{\lambda_m}), 1 \leq m \leq n, 1 \leq p(\lambda_m) \leq \infty, p(\lambda_m) = inf\{\alpha \in [1, \infty] : A_{\lambda_m} \in C\alpha H\lambda m$, $supp\lambda m: 1\leq m\leq n=\infty$ and for some $k\in N, A\lambda k\bar\in C\infty H\lambda k$, then for every $1\leq p < \infty$, $A \bar\in C_p(H)$.

**Proof.** The validity of the claims (1) and (2) are clear. Prove third assertion of theorem. If the operator $A \in C_p(H)$, then the series $\sum_{k=1}^\infty s_k^p(A)$ is convergent. In this case by the first proposition of this theorem and important theorem on the convergence of the rearrangement series it is obtained that the series $\sum_{m=1}^n \sum_{q=1}^\infty s_q^p(A_{\lambda_m})$ is convergent.

On the contrary, if the series $\sum_{m=1}^{n}\sum_{q=1}^{\infty} s_q^p (A_{\lambda_m})$ is convergent, then the series $\sum_{k=1}^{\infty} s_k^p (A)$ being a rearrangement of the above series, is also convergent. So $A \in C_p(H)$.

Now prove (4). If for every $m, 1 \leq m \leq n, \|A_{\lambda_m}\| \leq 1$, then from the inequality

$$\sum_{m=1}^{n}\sum_{q=1}^{\infty} s_q^p (A_{\lambda_m}) \leq \sum_{m=1}^{n}\sum_{q=1}^{\infty} s_q^{p(\lambda_m)} (A_{\lambda_m}) < \infty$$

and first claim the validity of this assertion is clear. Now consider the general case. In this case the operator $A$ can be written in form

$$A = CB, C = \oplus_{m=1}^{n}(1 + \|A_{\lambda_m}\|)E_m, B = \oplus_{m=1}^{n}\left(\frac{A_{\lambda_m}}{1+\|A_{\lambda_m}\|}\right).$$

Then $C \in B(H)$. On the other hand, since $\|B_m\| \leq 1, 1 \leq m \leq n$, and

$$\sum_{m=1}^{n}\sum_{q=1}^{\infty} s_q^{p(\lambda_m)}(B_m) = \sum_{m=1}^{n}\sum_{q=1}^{\infty} \frac{s_q^{p(\lambda_m)}(A_{\lambda_m})}{(1+\|A_{\lambda_m}\|)^{p(\lambda_m)}} \leq \sum_{m=1}^{n}\sum_{q=1}^{\infty} s_q^{p(\lambda_m)}((A_{\lambda_m}) < \infty,$$

then from the (3) of this theorem it implies that $B \in C_p(H)$ with $p = sup\{p(\lambda_m): 1 \leq m \leq n\}$.

Therefore, by the important theorem of the operator theory $A = CB \in C_p(H)$ [1].

Furthermore, by using proposition (2) of this theorem it is easy to prove the claim (5). On the other hand, the claim (6) is one of the corollary of (5).

**Remark 3.6.** Note that for the some $\lambda^* \in \Lambda^*$ in representation $\sigma_p(A) = \cup_{m=1}^{n} \sigma_p(A_{\lambda_m})$ in Theorem 3.4. may be hold

$$card\left[\sigma_p(A_{\lambda^*}) \cap \sigma_p(A)\right] < \infty$$

In these situations corresponding conditions for such index in the Theorem 3.5(3-6) may be omitted, for example, as in the following assertion.

**Theorem 3.7.** Assumed that $A = (A_\lambda), \lambda \in \Lambda, A = \int_\Lambda^\oplus A_\lambda \, d\mu(\lambda)$ in the Hilbert space $H=(H_\lambda), \lambda \in \Lambda$, $H = \int_\Lambda^\oplus H_\lambda \, d\mu(\lambda)$, $A_\lambda \in C_\infty(H_\lambda), \lambda \in \Lambda$ and $A \in C_\infty(H)$. In this case there exist countable subset $\Lambda^* = \{\lambda_1, \lambda_2, \lambda_3, \dots, \lambda_n\}, n \leq \infty$ of $\Lambda$ such that $\{s_k(A): k \geq 1\} = \cup_{m=1}^{n}\{s_q(A_{\lambda_m}): q \geq 1\}$

If

$card\{\lambda^* \in \Lambda^*: card\left[\sigma_p(A_{\lambda^*}) \cap \sigma_p(A)\right] < \infty\} < \infty$, $A_{\lambda_m} \in C_{p(\lambda_m)}(H_{\lambda_m}), 1 \leq p(\lambda_m) < \infty$, $\lambda_m \in \Lambda^*, \lambda_m \in \Lambda^{**} = \{\lambda^* \in \Lambda^*: card\left[\sigma_p(A_{\lambda^*}) \cap \sigma_p(A)\right] < \infty\}, 1 \leq p = sup\{p(\lambda_m): \lambda_m \in \Lambda^*\setminus\Lambda^{**}\} < \infty$

and $\sum_{\lambda_m \in \Lambda^*\setminus\Lambda^{**}}\left(\sum_{k=1}^{\infty} s_k^{p(\lambda_m)}(A_{\lambda_m})\right) < \infty$, then $A = (A_\lambda) \in C_p(H)$.

# 4. Power and Polynomially Boundednessity of the Direct Sum Operators

In this section let us $\Lambda = N, \Sigma = P(N)$ and $\mu$ is the counting measure. Here a connection of power (and polynomially) boundedness property of the direct sum operators in the direct sum Hilbert spaces and its coordinate operators were established. In advance, give some necessary definitions for the later.

**Definition 4.1**[25, 26]. Let $\mathcal{H}$ be any Hilbert space.

(1). An operator $T \in B(\mathcal{H})$ is called power bounded ($T \in PW(\mathcal{H})$), if there exist a constant $M(\geq 1)$

such that for any $n \in N$ is satisfied $\|T^n\| \leq M$. (3.1)

(2) Operator $T \in B(\mathcal{H})$ is called polynomially bounded ($T \in PB(\mathcal{H})$), if there exist a constant $M(\geq 1)$ such that for any polynomial $p(.)$ is satisfied $\|p(T)\| \leq M\|p\|_\infty$, (3.2) where $\|p\|_\infty = \sup\{|p(z)| : z \in \mathbb{C}, |z| \leq 1\}$.

(3) The smallest number $M$ satisfying (3.1) (resp.(3.2)) is called the power bound (resp. polynomial bound) of the operator $T$ and will be denoted by $M_w(T)$ (resp. $M_p(T)$).

Before of all note that the following theorem is true.

**Theorem 4.2.** If $A = \oplus_{n=1}^\infty A_n \in PW(H), H = \oplus_{n=1}^\infty H_n$, then for every $n \geq 1, A_n \in PW(H_n)$.

This result is a one of the corollary of following equality

$$\sup\nolimits_{m\geq 1}(\sup\nolimits_{n\geq 1}\|A_n^m\|) = \sup\nolimits_{n\geq 1}(\sup\nolimits_{m\geq 1}\|A_n^m\|) < \infty.$$

In general, the inverse of last assertion may be not true.

**Example 4.3.** Let $H_n = L^2(-1,1), \alpha_n \in \mathbb{R}, n \geq 1, \sup_{n\geq 1}|\alpha_n| = \infty$, $A_n : L^2(-1,1) \to L^2(-1,1)$,

$$A_n f(x) = \alpha_n \int_{-x}^x f(t)dt, \ n \geq 1, H = \oplus_{n=1}^\infty H_n, A = \oplus_{n=1}^\infty A_n.$$

It is easy to see that $\|A_n\| = 4|\alpha_n|/\pi$ and $A_n^2 = 0, n \geq 1$. Consequently,

$\sup_{m\geq 1}\|A^m\| = \sup_{m\geq 1}(\sup_{n\geq 1}\|A_n^m\|) = \sup_{n\geq 1}\|A_n\| = \sup_{n\geq 1} 4|\alpha_n|/\pi = \infty$.

Hence, for any $n \geq 1, A_n \in PW(H_n)$, but $A = \oplus_{n=1}^\infty A_n \bar\in PW(H)$.

**Example 4.4.** Let $H_n = \mathbb{C}^2$,

$$A_n = \begin{pmatrix} 0 & 0 \\ \alpha_n & 0 \end{pmatrix}, \alpha_n \in \mathbb{C}, n \geq 1, \sup_{n\geq 1}|\alpha_n| = \infty, A_n : \mathbb{C}^2 \to \mathbb{C}^2 \text{ and } H = \oplus_{n=1}^\infty H_n, A = \oplus_{n=1}^\infty A_n.$$

In this case $\|A_n\| = |\alpha_n|$ and $A_n^2 = 0, n \geq 1$, i.e. for every $n \geq 1, A_n \in PW(\mathbb{C}^2)$, but

$\sup_{m\geq 1}\|A^m\| = \sup_{m\geq 1}(\sup_{n\geq 1}\|A_n^m\|) = \sup_{n\geq 1}\|A_n\| = \sup_{n\geq 1}|\alpha_n| = \infty$.

Therefore, $A \bar\in PW(H)$.

Actually, it is true the following result.

**Theorem 4.5.** $A \in PW(H)$ if and only if $A_n \in PW(H_n), n \geq 1$ and $\sup_{n\geq 1} M_w(A_n) < \infty$.

**Proof.** If $A \in PW(H)$, then from the following relation

$$\sup_{m\geq 1}\|A^m\| = \sup_{m\geq 1}(\sup_{n\geq 1}\|A_n^m\| = \sup_{n\geq 1}(\sup_{m\geq 1}\|A_n^m\| < \infty$$

it is implies that $\sup_{m\geq 1}\|A_n^m\| < \infty$ for each $n \geq 1$. From this it is determined that for any $n \geq 1, A_n \in PW(H_n)$.

On the other hand, it is clear that for each $n \geq 1$

$$\sup_{m\geq 1}\|A_n^m\| \leq \sup_{m\geq 1}(\sup_{n\geq 1}\|A_n^m\| = \sup_{n\geq 1}(\sup_{m\geq 1}\|A_n^m\| = \sup_{m\geq 1}\|A^m\| \leq M_w(A) < \infty.$$

Therefore, $\sup_{n\geq 1} M_w(A_n) \leq M_w(A) < \infty$.

On the contrary, if for any $n \geq 1, A_n \in PW(H_n), \sup_{m\geq 1}\|A_n^m\| \leq M_w(A_n)$ and $\sup_{n\geq 1} M_w(A_n) < \infty$, then from the equality

$$\sup_{m\geq 1}\|A^m\| = \sup_{m\geq 1}(\sup_{n\geq 1}\|A_n^m\| = \sup_{n\geq 1}(\sup_{m\geq 1}\|A_n^m\| \leq \sup_{n\geq 1} M_w(A_n) < \infty$$

it is obtained that $A \in PW(H)$.

Now polynomially boundness property of the direct sum operators will be investigated. In advance, note that the following proposition is true.

**Theorem 4.6.** If $H = \oplus_{n=1}^{\infty} H_n$, $A = \oplus_{n=1}^{\infty} A_n$ and $A \in PB(H)$, then for every $n \geq 1, A_n \in PB(H_n)$.

Unfortunately, the inverse of last theorem may be not true in general.

**Example 4.7.** Let $H_n = L^2(-1,1), \alpha_n \in \mathbb{R}, n \geq 1, \alpha_n \geq \pi/4, \sup_{n\geq 1}|\alpha_n| = \infty$, $A_n: L^2(-1,1) \to L^2(-1,1)$,

$$A_n f(x) = \alpha_n \int_{-x}^{x} f(t) dt, n \geq 1,$$

It is known that $A_n \in C_{\infty}(H_n)$ and $A_n$ is a nilpotent operator with power of nilpotency 2, for any $n \geq 1$. In this case for any polynomial function $p(z) = \sum_{k=0}^{q} a_k z^k, z \in \mathbb{C}, q = 0,1,2, \ldots$ we have

$$\|P(A_n\| \leq |a_0| + |a_1|\|A_n\| \leq 4|\alpha_n|/\pi (|a_0| + |a_1|) = \|A_n\|\|p\|_{\infty}, n \geq 1.$$

In the other words, for every $n \geq 1, A_n \in PB(H_n)$. Unfortunately, for the polynomial $p_*(z) = z, z \in \mathbb{C}$ we have

$$\|p_*(A)\| = \sup_{n\geq 1}\|p_*(A_n)\| = \sup_{n\geq 1}\|A_n\| = \sup_{n\geq 1} 4|\alpha_n|/\pi = \infty,$$

i.e. $A \bar{\in} PB(H)$.

But in general case the following result is true.

**Theorem 4.8.** Let $H = \oplus_{n=1}^{\infty} H_n$, $A = \oplus_{n=1}^{\infty} A_n$ and $A \in B(H)$. In order to $A \in PB(H)$ the necessary and sufficient conditions are $A_n \in PB(H_n)$ for every $n \geq 1$ and $\sup_{n\geq 1} M_p(A_n) < \infty$.

**Proof.** Assumed that for every polynomial p(.) $\|p(A_n)\| \leq M_p(A_n)\|p\|_{\infty}$. and $\sup_{n\geq 1} M_p(A_n) < \infty$. In this case since for every polynomial function $p(.)$ $P(A) = P(\oplus_{n=1}^{\infty} A_n) = \oplus_{n=1}^{\infty} P(A_n)$, then $\|p(A_n)\| = \sup_{n\geq 1}\|P(A_n)\|$. From last relation it is obtained that $\|p(A)\| \leq \sup_{n\geq 1} M_p(A_n)\|p\|_{\infty}$. Hence $A \in PB(H)$.

Now let us $A \in PB(H)$, i.e. for any $n \geq 1$ and polynomial p(.) it is valid

$$\| p(A_n) \| \leq \| p(A) \| = \sup_{n \geq 1} \| P(A_n)\| \leq \sup_{n \geq 1} M_p(A_n) \|p\|_\infty.$$

Then it is clear that $A_n \in PB(H_n), n \geq 1$. In the other hand, from last equality it implies that for every $n \geq 1$ is hold $M_p(A_n) \leq M_p(A) < \infty$. Hence, $\sup_{n \geq 1} M_p(A_n) < \infty$. This completes the proof of the theorem.

**Acknowledgment**

The authors are grateful to G.İsmailov(student of Trabzon Kanuni Anadolu High School) for his helping suggestion to english version and other technical discussion.